\documentclass[11pt]{amsart}
\usepackage{euler,amsmath,amssymb,amscd}
\usepackage{epsfig}
\usepackage{array, multirow}
\usepackage[usenames]{color}

\newtheorem{theorem}{Theorem}

\newtheorem{proposition}{Proposition}
\newtheorem{prop}{Proposition}

\newtheorem{lemma}{Lemma}

\theoremstyle{definition}

\newtheorem{remark}{Remark}

\newtheorem{ex}{Exercise}
\newtheorem{definition}{Definition}

\def\mbb{\mathbb}

\def\mcl{\mathcal}

\def\ten{\otimes}

\def\tu{\textup}

\def\lan{\langle}
\def\ran{\rangle}

\def\a{\alpha}

\def\d{\delta}
\def\D{\Delta}
\def\e{\epsilon}

\def\D{\Delta}

\def\bP{\mbb P}

\def\bG{\mathbb G}

\def\inj{\hookrightarrow}

\def\spec{\tu{Spec\,}}

\def\Proj{\tu{\textbf{Proj}\,}}

\def\om2{\omega^{\ten 2}}

\def\Mg{\overline{M}_g}

\def\FMps{\overline{{\bf \mcl M}}^{ps}_g}
\def\FMg{\overline{{\bf \mcl M}}_g}

\def\inj{\hookrightarrow}
\def\GL{GL}

\def\Pic{Pic}

\def\Sym{Sym}

\def\M{\bar{M}}
\def\cO{\mathcal O}
\def\Mps{\overline{M}_3^{ps}}
\def\dps{\delta^{ps}}

\def\cM{\mathcal{M}}

\def\ra{\rightarrow}

\def\bar{\overline}

\def\codim{\textup{codim}\,}

\input xy
\xyoption{all}

\input epsf
\def\bZ{\mathbb Z}

\def\bQ{\mathbb Q}
\DeclareMathOperator{\Hilb}{Hilb}

\DeclareMathOperator{\Chow}{Chow}
\DeclareMathOperator{\Aut}{Aut}
\DeclareMathOperator{\Stab}{Stab}
\DeclareMathOperator{\Spec}{Spec}
\DeclareMathOperator{\Sing}{Sing}
\DeclareMathOperator{\Cl}{Cl}
\DeclareMathOperator{\Def}{Def}

\begin{document}

\title[GIT constructions of log canonical models]{GIT constructions of log canonical models of $\Mg$}
\date{\today}

\author{Jarod Alper}
\address{
Departamento de Matem\'aticas\\
Universidad de los Andes\\
Cra 1 No. 18A-10\\
Edificio H\\
Bogot\'a, 111711, Colombia\\
} \email{jarod@uniandes.edu.co}

\author{Donghoon Hyeon}
\address{
Department of Mathematics\\Postech\\ Pohang, Gyungbuk
790-784\\ Republic of Korea
} \email{dhyeon@postech.ac.kr}

\maketitle

\section{Introduction}
The purpose of this article is to give an overview of the construction of compact moduli spaces of curves from the viewpoint of the log minimal model program for $\Mg$ (now coined as the ``Hassett-Keel" program).  We will provide an update on new developments and discuss further problems. We have attempted to complement recent articles by Fedorchuk and Smyth \cite{FS} and Morrison \cite{Mor}, and as a result our focus is put on the GIT construction of moduli spaces  using {\it low degree Hilbert points} of curves. This method is expected to produce new compact moduli spaces of curves with increasingly worse singularities, and recent work by Ian Morrison and Dave Swinarski \cite{MS} is one solid step forward in this direction.  The low-degree Hilbert quotients conjecturally realize various log canonical models of $\Mg$
\[
\Mg(\a) = Proj \oplus_{m\ge 0} \Gamma(\Mg, m(K_{\FMg} + \a \d))
\]
as moduli spaces, for certain values of $\a \in [0,1].$ 
In the excellent survey \cite{FS}, Fedorchuk and Smyth place the Hassett-Keel program in the larger context of birational geometry and classification of modular compactifications.  On one hand, the Hassett-Keel program has stirred much excitement in studying and constructing various new compactifications; in fact, the ambitious paper \cite{Smy} provided a complete classification of stable modular compactifications.  On the other hand, the log canonical models of $\M_3$ essentially account for all known compactifications for $M_3$, and gives a complete Mori decomposition of the restricted effective cone; see \cite{HL2}. These entwined views are nicely discussed in \cite{FS}.  The GIT construction of moduli spaces of curves with applications toward the Hassett-Keel program is carefully reviewed in Morrison's survey \cite{Mor} which we will be referencing frequently.

Here we will be mainly interested in the sort of moduli problems that are ``mildly non-separated" i.e. the ones involving moduli functors that are not separated but admit a projective moduli space nonetheless.  In the language of \cite{ASvdW}, these moduli problems are referred to as \emph{weakly separated}; there may be many ways to fill in a family of curves over a punctured disc but there is a unique limit which is closed in the moduli stack. 
GIT quotient stacks of the form $[X^{ss}/G]$ are the model examples. Moduli spaces of sheaves or more generally decorated sheaves as GIT quotients generally have strictly semistable points (and are therefore non-separated) but are weakly separated.  The study of moduli of curves has  revolved around the beautiful compactification $\Mg$ of Deligne-Mumford stable curves and weakly separated functors have become prominent only recently coinciding with their appearance in the Hassett-Keel program  for $\Mg$ \cite{HL1, HM, HH2}. Likewise, although GIT has been around for decades, the construction of compact moduli spaces of curves in which non-isomorphic curves are identified has become a topic of pronounced interest only recently.\footnote{It is worth mentioning here  that Compact moduli spaces of pairs $(C, E)$ consisting of a vector bundle $E$ over a curve $C$  generally have strictly semistable points \cite{GM, Caporaso,  Pandharipande}}
  We will elaborate on how GIT can be carried out in such problems. Once we have chosen a suitable parameter space, the construction proceeds roughly as follows: 
\begin{enumerate}
\item Establish the GIT stability of nonsingular curves;

\item Destabilize undesirable curves such as: Degenerate and non-reduced curves, badly singular curves and special subcurves such as tails and chains;

\item Enumerate all potentially semistable curves with infinite automorphisms, and compute their basin of attraction and possible semistable replacements;

\item To finalize the construction, prove that potentially semistable curves are semistable. The argument heavily depends on the semistable replacement theorem and the existing moduli spaces such as $\Mg$ and $\Mg^{ps}$ \cite{Sch}. Basically, given a potentially semistable curve $C$, the proof entails arguing that all semistable replacements of the Deligne-Mumford (or pseudo-stable) stabilization of $C$ is in a basin of attraction of a curve $C_0$ that is strictly semistable with respect to a one-parameter subgroup coming from $\Aut(C_0)$.
\end{enumerate}
This strategy is taken in the constructions in \cite{HL1}, \cite{HM} and \cite{HH2}.  Establishing the stability of smooth curves is definitely the most difficult step, and it was accomplished in \cite{Mum, Gie} for curves and in \cite{Swin} for weighted pointed curves. 
The instability analysis of degenerate curves, non-reduced curves and badly singular curves has become rather standard following the work of Mumford and Gieseker. Although the configuration of  destabilizing special subcurves in \cite{HH2} are more intricate and their instability analysis is accordingly more involved, destabilizing tails and bridges appear in the work of Gieseker, Mumford and Schubert \cite{Gie, Mum, Sch}: stable curves do not admit rational tail and rational bridges, and furthermore pseudo-stable curves do not admit elliptic tails; they are replaced by an ordinary cusp.  But as we will see, these standard steps fail when applied to low-degree Hilbert stability, and addressing these issues will be a major focus of this article. 
The third step above is actually a new feature of \cite{HL1, HM, HH2} since $\Mg$ and $\Mg^{ps}$ have only curves with finite automorphisms.  In the moduli spaces of curves $\Mg^{cs}$ and $\Mg^{hs}$ where curves with at worst nodes, cusps and tacnodes are parameterized, there are certain maximally degenerate curves (with infinite automorphisms) corresponding to minimal orbits in the Hilbert or Chow space; all other curves admits a unique isotrivial specialization to a maximally degenerate curve.  This is a very useful structural property that is central in the construction of a projective moduli space of weakly separated moduli functors.

Much of the material in this article has been drawn from other papers including the two nice surveys \cite{FS, Mor}. 
\begin{itemize}
\item In Section~\ref{S:parameter}, we explain how a parameter space can be chosen for the GIT setup where the aim is to construct a log canonical model $\Mg(\a)$ of curves with prescribed singularities. We follow expositions in \cite{H} and \cite[Section~2.4] {FS}.

 \item  In Section~\ref{S:finite}, we will briefly review Kempf's theory \cite{Kempf} and the Gr\"obner techniques developed in \cite{MS}. Using them, we provide two arguments to show that a tri-canonical  genus two rational curve with two ordinary cusps a has semistable $2$nd Hilbert point.  This curve played a central role in the stability analysis of tri-canonical genus two curves \cite{HL1}. These computations are new. 
 
 \item In Section~\ref{S:degenerate}, by using Gr\"obner basis technique \cite{BM, HHL}, we give a short proof of instability of $m$th Hilbert points of reduced degenerate curves, for any degree $m$. In the non-reduced case, the same proof  also gives instability for sufficiently large $m$ as well as an effective lower bound for $m$.   Non-reduced curves are in fact expected to appear in log canonical models arising from the GIT of  low degree Hilbert points. 
 
\item Although \cite{Gie} and \cite{Mum} provide  fundamental ideas on picking destabilizing one-parameter subgroups, the dimension estimation method in \cite{Gie} is ill suited for computing Hilbert-Mumford index of low degree Hilbert points since we do not have the vanishing of higher cohomology.
Instead, we let the undesirable curve specialize to a degenerate configuration and carry out the exact Hilbert-Mumford index computation there. This tells us precisely at which $m$ the $m$th Hilbert point is semistable.   In Section~\ref{S:special}, we first explain how to figure out which curves specialize to a given maximally degenerate configuration curve with $\bG_m$ action (Section~\ref{S:basin}).   We then apply these ideas to provide two methods to give precise predictions for  which $m$ an $A_{2k+1}$-singularity would be $m$-Hilbert (semi)stable. 

\item In Section~\ref{S:factorial}, we provide an \'etale local description of the flip  $\Mg(9/11) \to \Mg(7/10) \leftarrow \Mg(7/10-\e)$ introduced in \cite{HH2} based on the ideas of \cite{ASvdW}.  Moreover, we show that $\Mg(7/10-\e)$ are not $\bQ$-factorial.  However, by scaling a generic boundary divisor rather than the ``democratic'' boundary, one expects $\bQ$-factorial flips.  For a generic scaling, we provide a conjectural description of the critical values and the corresponding base loci.
\end{itemize}

\noindent We work over an algebraically closed field $k$ of characteristic zero. 

\subsection{ Acknowledgments}
We would like to thank Dave Swinarski for helpful discussions on Section~\ref{S:finite} and for carrying out the state polytope computation with Macaulay 2. He also read a preliminary version carefully and gave us many useful comments.
  We also thank Maksym Fedorchuk for suggesting the Hilbert-Mumford index analysis at the end of Section~\ref{S:finite}, and David Smyth for helpful discussions.

\section{Parameter spaces} \label{S:parameter}
Hilbert schemes and Chow varieties are natural candidates for parameter spaces in the construction of moduli spaces of curves. While Chow varieties come with a canonical linearization, Hilbert schemes admit a family of  linearizations, and one should choose the linearization depending on the desired properties of the space that one wishes to construct.  From the point of view of the Hassett-Keel program, we would like to construct a moduli space on which $13 \lambda - (2-\alpha)\delta$ is ample, for a given $\alpha \in [0,1]\cap \bQ$. 

Given a nonsingular curve $C$ of genus $g \ge 2$,  $V:=H^0(\omega_C^{\otimes n})$ is a $k$-vector space of dimension $N+1 = (2n-1)(g-1)$ if $n \ge 2$ or $g$ if $n=1$.  The induced embedding $
C \hookrightarrow \bP(V)
$
defines the Chow point $\Chow(C)$ and the Hilbert point $[C]$. Taking the closure of the locus of such points (corresponding to nonsingular curves) in the Chow variety and the Hilbert scheme, we obtain our parameter spaces $\Chow_{g,n}$ and $\Hilb_{g,n}$, respectively. 

Recall that $\Hilb_{g,n}$ admits a family of embeddings
\[
\phi_m: \Hilb_{g,n} \hookrightarrow Gr(P(m), S^mV^*) \hookrightarrow \bP\left(\bigwedge^{P(m)}S^mV^*\right)
\] 
when $m$ is sufficiently large ($\ge$ Gotzmann number \cite{Gotzmann}). 
\begin{definition} The image $[C]_m$ of $[C]\in \Hilb_{g,n}$ under $\phi_m$ is called the $m$th Hilbert point of $C$, and $C$ is said to be $m$-Hilbert stable (resp., semistable, unstable) if $[C]_m$ is stable (resp., semistable, unstable) with respect to the natural $SL(V)$ action linearized by the embedding $\phi_m$. $C$ is Hilbert stable (resp., semistable, unstable), or asymptotically Hilbert stable (resp., semistable, unstable) if it is $m$-Hilbert stable (resp., semistable, unstable) for $m\gg 0$. 
\end{definition}
If we are concerned with the $m$-Hilbert stability for fixed $m$ smaller than the Gotzmann number so that $\phi_m$ is not defined on the whole of $\Hilb_{g,n}$, the GIT problem will be referred to as {\it finite Hilbert stability} problem. 

\

 According to Mumford's Grothendieck-Riemann-Roch computation \cite[Theorems~5.10, 5.15]{Mum}, the canonical polarization on $\Chow_{g,n}$ is positive rationally proportional to 
 \[
\begin{cases}
(4g+2)\lambda - \frac g2 \delta, \quad n = 1 \\(6n-2)\lambda -\frac n2 \delta, \quad \mbox{otherwise}
\end{cases}
\]  
where $\lambda$ is the determinant of the Hodge bundle and $\delta$ is the divisor of the singular curves. 
Mumford's formula, combined with a basic Chern class computation, also shows that  $\Lambda_{m,n}:=\phi_m^*(\cO(+1))$ on $\Hilb_{g,n}$ is a positive rational multiple of
\begin{equation}\label{E:Lmn}
\begin{cases}
\lambda + (m-1) [((4g+2)m-g+1)\lambda - \frac{gm}2 \delta], \quad n=1\\
(6mn^2-2mn-2n+1)\lambda -\frac{mn^2}2 \delta, \quad n>1
\end{cases}
\end{equation}
\cite[P28. Equation~(5.3)]{HH2}. From this, we see that $\Lambda_{m,n}$ is an $SL(V)$-linearized $\bQ$-Cartier divisor. This allows us to define $m$-Hilbert (semi)stability for rational $m$: 
\begin{definition} Given a positive rational number $m$, we say that $[C] \in \Hilb_{g,n}$ is $m$-Hilbert (semi,un)stable if it is GIT (semi,un)stable with respect to $\ell\Lambda_{m,n}$ for sufficiently large $\ell \in \bZ$.
\end{definition}

Equation $(1)$ suggests that the log canonical model $\Mg(\a)$ may be obtained from the GIT of $\Hilb_{g,n}$ with a suitable linearization $\Lambda_{m,n}$.  More precisely, suppose we know that:
\begin{itemize}
\item $\Hilb_{g,n}^{ss,m} \neq \emptyset$, where $\Hilb_{g,n}^{ss,m}$ denotes the semistable locus with respect to the linearization $\Lambda_{m,n}$.
\item The locus  in $\Hilb_{g,n}^{ss,m}$ of curves which are not Deligne-Mumford stable has codimension at least 2.
\end{itemize}
Then the map $\overline{M}_g \dashrightarrow \Hilb_{g,n}/\!\!/SL(V)$ is a birational contraction.  The natural line bundle on $\Hilb_{g,n}/\!\!/SL(V)$ descended from $\Lambda_{m,n}$ pulls back to a linear combination of $\lambda$ and $\delta$ on $\overline{M}_g$ as specified in Equation $(1)$ modulo exceptional divisors.  It follows that $\Hilb_{g,n}/\!\!/SL(V)$ is identified with $\Mg(\a)$ where $\a = 2 - 13/s$ and $s$ is the slope of the divisor in Equation $(1)$ since the addition of positive multiples of exceptional divisors has no effect on the section ring. 

For instance, when $n=2$, $\Lambda_{m,n}$ is proportional to  $13\lambda - (2 - \a)\delta$.  For meaningful values of $m$ and $\a$,  there is a one-to-one order preserving correspondence:
\begin{equation}\label{E:anm}
m(\a) = \frac{3(2-\a)}{2(7-10\a)}, \quad \a(m) = \frac{14m-6}{20m-3}.
\end{equation}
 By computing the $m(\a)$-Hilbert semistable curves, we can predict which singular curves may appear in $\Mg(\a)$. 
\cite{FS} also explains how one can predict the singularities that appear in $\Mg(\a)$ as constructed by GIT, by using character theory. 
Table~\ref{T:1} (reproduced from \cite{FS} at the end of this article) nicely summarizes what has been shown and what we expect. The last four lines of the table are conjectural with a few partial results supporting them. In \cite{HL4} is given a construction of a birational map $\M_4^{hs} \to \M_4(2/3)$ that contracts the locus of Weierstrass genus $2$ tails, which is precisely the variety of stable limits of $A_4$-singularities. This is consistent with the prediction that Weierstrass genus two tails would be replaced by $A_4$-singularities in $\Mg(2/3)$. In \cite{ASvdW}, Alper, Smyth and van der Wyck takes a GIT-free approach in constructing the log canonical models. Their plan is to first define the moduli functor with desired properties, show that it is deformation open and weakly proper, and  prove that a projective coarse moduli space exists. They accomplish the first step of this ambitious program in \cite{ASvdW}.

\section{ Finite Hilbert Stability}\label{S:finite}
First, let us recall the classic result on stability of smooth curves.
\begin{theorem} \cite{Gie, Mum}
 A smooth curve of genus $g$ embedded by a complete linear system of degree $d \ge 2g+1$ is both Hilbert stable and Chow stable. 
 \end{theorem}
The proof of this theorem is thoroughly reviewed in \cite{Mor} and we will not discuss it here. Instead, 
we will address the question of finite Hilbert stability;  see Table~\ref{T:1} for the correlation between the integer $m$ and the log canonical models $\bar{M}_g(\alpha)$.   We emphasize that finite Hilbert stability is completely different from the asymptotic case  as we do not have the vanishing of higher cohomology and the weight space dimension estimation method is not feasible. 
We will explain how one can prove the finite Hilbert semistability of certain singular curves whose automorphism groups satisfy a multiplicity  condition. 
 Establishing the semistability of any smoothable curve of course  implies the semistability of a generic smooth curve, and though  it is certainly far short of the stability of {\it all} smooth curves, 
 we allow ourselves to be content with it for now. To complete the stability proof, we would have to resort to a semistable replacement argument.

Recall the setup: the parameter space is the closure in the Grassmannian of the image of 
\[
Hilb_{g,2}^\circ \inj Gr(P(m), S^mV^*)
\]
where the superscript $\circ$ denotes the locus of $m$-regular curves. Immediately we realize that even destabilizing unwanted objects is perhaps more subtle as the boundary points may not even correspond to curves. 

\subsection{Semistability of the bicuspidal rational curve.}
\begin{figure}[ht]
\centerline{\scalebox{0.45}{\psfig{figure=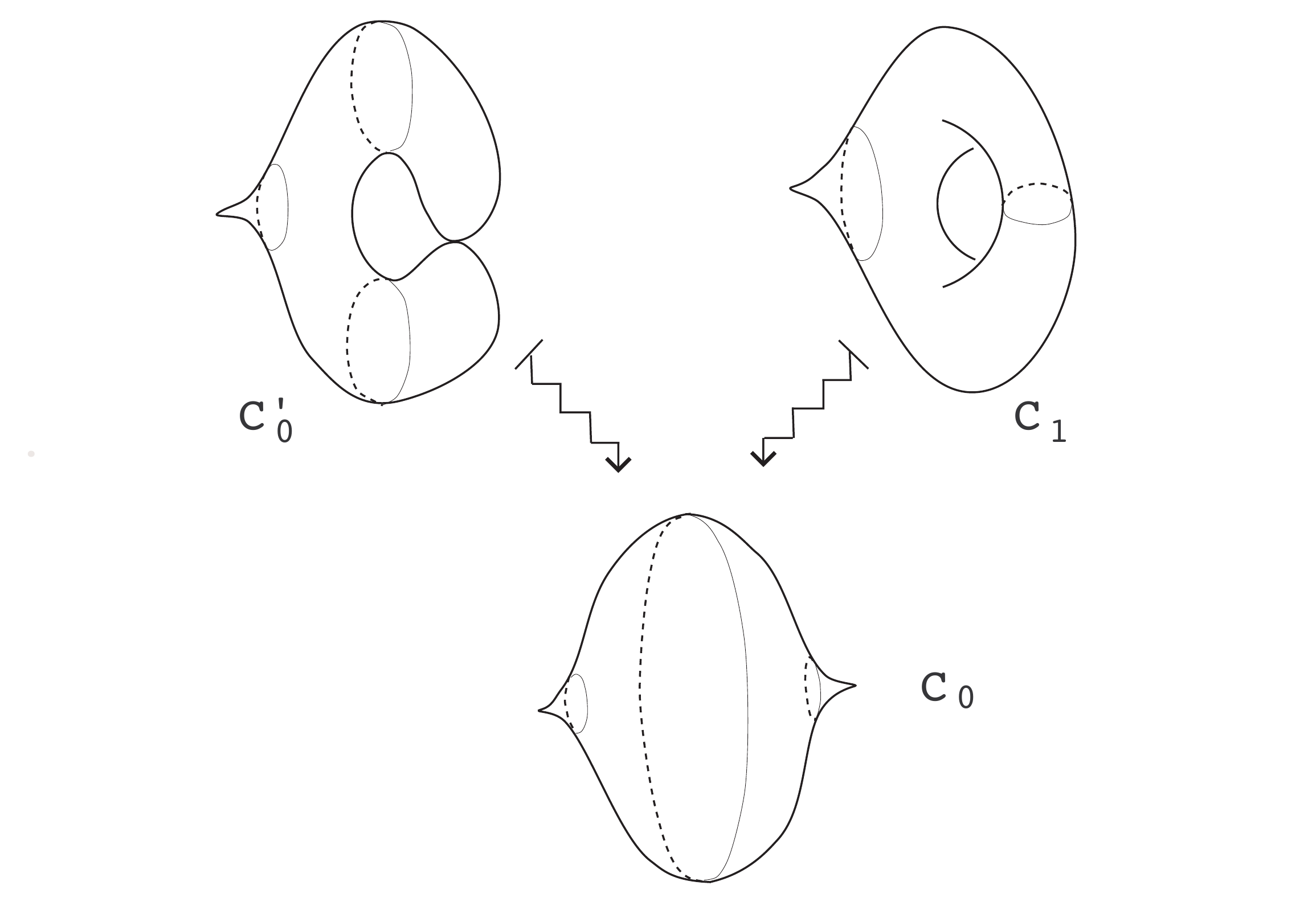}}}
\caption{Basin of attraction of a bicuspidal genus two curve}
\label{F:bicuspidal}
\end{figure}

In this section, we will sketch the technique from \cite{MS} and use it to show that a tri-canonical genus two rational curve $C_0$
 with two ordinary cusps is $2$-Hilbert  semistable. 
\begin{remark} Definition~4.8 in \cite{MS} assumes that the automorphism group of the variety in question has a finite subgroup acting by linear automorphisms, but their method works without the finiteness assumption as we shall demonstrate below.
\end{remark}
 Note that $C_0$ is unique up to isomorphism. 
This particular curve played an important role in the construction of $\M_2^{ps}$ in \cite{HL1}. It is the unique strictly semistable pseudo-stable curve of genus two with infinite automorphisms. Other strictly semistable curves specialize to it via the $\bG_m$ action coming from $\Aut(C_0)$ (Figure~\ref{F:bicuspidal}). 
The tri-canonical image  of $C_0$ may be parametrized by
\[
\begin{array}{cccc}
\nu : \bP^1 & \to & \bP^4 =: \bP(V)  \\
\left[s,t\right] & \mapsto & [s^6, s^4t^2, s^3t^3, s^2t^4, t^6]
\end{array}
\]
from which we compute its defining ideal:
\[
\lan x_3^2-x_1 x_4, x_1 x_3-x_0 x_4, x_2^2-x_0 x_4, x_1^2-x_0
x_3\ran.
\]
For $\a \in \bG_m$, the automorphism $[s, t] \mapsto [\a s, t]$ of $\bP^1$ induces an automorphism $\phi_\a \in \GL(V)$ of $C_0$, and we denote by $\Gamma$ the subgroup of $\GL(V)$ consisting of $\phi_a$'s. Note that $\Gamma$ is (isomorphic to) the identity component of the automorphism group of $C_0$. We have 
\[
\Gamma\simeq \bG_m \subset \Stab_{\GL(V)}([C_0]_m)
\]
and $\Gamma$ acts on $V$ with weights $(6,4,3,2,0)$ inducing the weight space decomposition
$
V = V_6 \oplus V_4 \oplus V_3 \oplus V_2 \oplus V_1
$
such that $\a.v = \a^m v$, $\forall (\a, v) \in \bG_m\times V_m$. 
The key point is that, since the weights are all distinct (this corresponds to the multiplicity free condition  \cite[Definition~4.5]{MS}), the  decomposition above
determines a unique maximal torus $T_\Gamma \subset \GL(V)$.  Note that the coordinates we use above are compatible with $T_\Gamma$. 

Now,  Kempf's theory  \cite[Theorem 3.4 and Corollary 3.5]{Kempf} says the following. 
\begin{enumerate}
\item 
If $[C_0]_m \in \Hilb_{2,3}$ were unstable, then there would be a {\it worst} one-parameter subgroup $\rho^\star$ whose associated parabolic subgroup $P$ compatible with the $\rho^\star$-weight filtration contains $\Stab_{\GL(V)}([C_0]_m)$. 

\item Also, if $T'$ is a maximal torus contained in $P$, then there is a 1-PS $\rho'$ of $T'$ such that 
\[
\mu([C_0]_m, \rho^\star) = \mu([C_0]_m, \rho').
\]
\end{enumerate}
But  $\Gamma \subset \Stab_{\GL(V)}([C_0]_m) \subset P$, $\Gamma$ preserves the flag of $P$ and by the complete reducibility of $\Gamma$, each step of the flag is a direct sum of $\Gamma$-weight spaces. Hence $T_\Gamma$ also preserves the flag of $P$ i.e. $T_\Gamma$ is a maximal torus of $\GL(V)$ contained in $P$, and Kempf's theory implies that there exists a 1-PS $\rho$ of $T_\Gamma$  such that 
$
\mu([C_0]_m, \rho^\star) = \mu([C_0]_m, \rho).
$

In summary, if $[C_0]_m$ were unstable, then it is destabilized by a one-parameter subgroup $\rho$ contained in our favored maximal torus $T_\Gamma$ compatible with the $\Gamma$-weight space decomposition. In other words, we only need to check the stability with a fixed basis compatible with $T_\Gamma$, for instance $(x_0, \dots, x_4)$ from before,
and we will do just that by employing
the systematic state polytope trick developed in  \cite{BM} (implemented by Dave Swinarski into the Macaulay 2 package \verb"statePolytope"). In the code below, \verb"statePolytope(2,I)" computes the 2nd state polytope of the ideal which has 10 vertices. The command \verb"isStable(2, I)" confirms that the $2$nd state polytope of $I$ contains the barycenter $(8/5, \dots, 8/5)$, so $C_0$ is $2$-Hilbert semistable.

\footnotesize
\begin{verbatim}
i1 : loadPackage("StatePolytope");

i2 : R=QQ[a,b,c,d,e];

i3 : I = ideal(-b*e+d^2,-a*e+b*d,-a*e+c^2,-a*d+b^2);

i4 : statePolytope(2,I)

o4 = {{1, 3, 0, 3, 1}, {1, 4, 0, 1, 2}, {2, 1, 0, 4, 1}, {2, 2, 0, 2, 2}, {1,
     --------------------------------------------------------------------------
     1, 2, 4, 0}, {2, 0, 2, 3, 1}, {2, 1, 2, 1, 2}, {0, 3, 2, 3, 0}, {0, 4, 2,
     --------------------------------------------------------------------------
     1, 1}, {1, 3, 2, 0, 2}}

o4 : List

i5 : isStable(2,I)

o5 = true
\end{verbatim}
\normalsize

In fact, the ideal is simple enough that the computation can be done easily by hand in the degree two case. Given a one-parameter subgroup $\bG_m \to SL_{4}$ with weights $(r_0, r_1, \dots, r_4)$, we shall prove that there exists a monomial basis of $H^0(C, \mathcal{O}_C(2))$ whose weights have non-positive sum.  One easily checks that the set $B_1$ (resp.,  $B_2$, $B_3$) consisting of all monomials except for $\{x_1^2, x_2^2, x_3^2, x_0 x_4\}$ (resp., $\{x_1x_4, x_0x_4, x_0x_3, x_1x_3\}$, $\{x_1x_4, x_0x_4, x_0x_3, x_2^2\}$) is a monomial basis of $H^0(C, \mathcal{O}_C(2))$.  If the sum of the weights of the basis $B_1$ is positive, then $r_0 + 2r_1+2r_2+2r_3+r_4 < 0$ or $r_1+r_2+r_3 < 0$.   If the sum of the weights of the basis $B_2$ is positive, then $2r_0+2r_1+2r_3+2r_4 < 0$ or $r_2 > 0$.   If the sum of the weights of the basis $B_3$ is positive, then $2r_0+r_1+2r_2+r_3+ 2r_4 < 0$ or $r_1 + r_3 > 0$.  Obviously, all three inequalities cannot hold at once, and  one of the bases $B_1, B_2$ or $B_3$ must have a non-positive sum of weights and we conclude that $C$ is $2$-Hilbert semistable.

\section{Unstable curves: Degenerate and non-reduced curves}\label{S:degenerate}
Mumford's numerical criterion can be used to destabilize an unwanted variety effectively. The following formulation of the Hilbert-Mumford index is given in \cite{BM}. It was heavily used in the work of Hassett, Hyeon, Lee and Morrison. A Macaulay 2 implementation is given in \cite{HHL}.

\begin{prop}
\label{prop:stabcrit2} Let $X \subset \bP(V)$ be a projective variety defined by a homogeneous ideal $I$.
The Hilbert-Mumford index of $[X]_m$ with
respect to a one-parameter subgroup $\rho : \bG_m \to \GL(V)$
with weights $r_0, r_1, \dots, r_N$ is given by
\begin{equation}\label{eqn:stabcrit}
\mu([X]_m, \rho) = \frac{m P(m)}{N+1} \sum r_i -\sum_{j=1}^{P(m)} wt_{\rho}(x^{a(j)})
\end{equation}
 where
$a(1),\ldots,a(P(m))$ index the monomials of degree $m$ not
contained in the initial ideal $in_{\prec_\rho}(I)$.
 In particular, $[X]_m\in \bP(\bigwedge^{P(m)} \Sym^mV)$ is
stable (resp. semistable) under the natural $\GL(V)$-action if and
only if for any one-parameter subgroup $\rho$ we have
$$\sum_{j=1}^{P(m)} wt_{\rho}(x^{a(j)})
< \frac{m P(m)}{N+1} \sum r_i \qquad (\text{resp.} \leq).$$
Here, $\prec_\rho$ denotes any fixed graded total order that refines the graded $\rho$-weight order. 
\end{prop}

Using this algorithm, we give the following short proof of the instability of degenerate varieties.  It implies that a reduced degenerate variety is $m$-Hilbert unstable for all $m$ whereas the results in the literature are of asymptotic nature.

\begin{lemma}  If $X \subset \bP^N$ is contained in the $r$th thickening of a hyperplane, then it is Hilbert unstable for $m > (N+1)(r-1)$.  In particular, if $X$ is a non-reduced degenerate variety, then $X$ is $m$-Hilbert unstable for all $m$. 
\end{lemma}

\begin{proof} Let $I$ be the homogeneous ideal of $X$.  The hypothesis implies that we may choose coordinates so that $x_0^r \in I$.
Let $\rho $ be the 1-PS with weights $(0,1,\dots,1)$. 
If $x_0^r$ divides a monomial $x^a$, then $x^a \in I$ and hence $x^a \in in_{\prec_\rho}(I)$. So any degree $m$ monomial not in the initial ideal has $\rho$-weight at least $m-r+1$. It follows that
$$\displaystyle{
\begin{array}{ccl}
\mu([X]_m, \rho) & = & - \sum_{\substack{x^a \not\in in_{\prec_\rho}(I)_m \\ ||a|| = m}} wt_\rho(x^{a}) + \frac{m P(m)}{N+1} \sum_{i=1}^{N+1} r_i \\
& \leq & - (m-r+1) P(m) + m P(m) \frac{N}{N+1}  \\
& = & P(m) \left( -\frac1{N+1} m + r - 1\right).
\end{array}}
$$
Here, $r_i$ are the weights of $\rho$ which sum up to $N$. It is evident that $\mu([X]_m, \rho) < 0$ for $m > (N+1)(r-1)$.
\end{proof}


\begin{remark} The ``standard" instability proof for non-reduced curves fails completely when applied to the finite Hilbert stability case, which it should, since non-reduced curves are expected to appear in $\Mg(5/9)$. This corresponds to the GIT quotient moduli space of bicanonical $3/2$-Hilbert semistable curves.
\end{remark}

\section{Unstable curves: Badly singular curves and special subcurves}\label{S:special}

From the work of Mumford and Gieseker, 
we have fairly standard if not fully systematic procedures to produce one-parameter subgroups for destabilizing curves with undesirable singularities.
But the proof in \cite{Gie} is ill suited for dealing with low degree Hilbert points since it is designed to work for general curves (with the singularity in question) and depends on the vanishing of higher cohomology. Instead, the proofs in \cite{HL1, HM, HH2} employ the following strategy. First we need to define the {\it basin of attraction}:

\begin{definition}\label{D:basin}
Let $X$ be a variety on which $\bG_m$ acts via $\rho:\bG_m \ra Aut(X)$
with fixed points $X^{\rho}$.
For each $x^{\star} \in X^{\rho}$, the {\em basin of attraction}
is defined
\[
A_\rho(x^\star) := \left\{ x\in X \, | \, \lim_{t \to 0} \rho(t).x = x^\star
\right\}.
\]
\end{definition}
Suppose that we are carrying out GIT : a reductive group $G$ acts on a (projective) variety $X$ linearly and $\rho : \bG_m \to G$ is a 1-PS fixing $x^\star$. If $x \in A_\rho(x^\star)$, then by definition $\mu(x,\rho) = \mu(x^\star, \rho)$. In particular, $x$ is unstable with respect to $\rho$ if and only if $x^\star$ is.  Moreover, if $x^\star$ is strictly semistable with respect to $\rho$, then $x$ is semistable if and only if $x^\star$ is;  see \cite[Lemma 4.3]{HH2}. 

The GIT strategy pursued in \cite{HL1, HM, HH2} for destabilizing a curve $C$ proceeds as follows: 
\begin{enumerate}
\item Find a degenerate configuration $C_0$ of $C$ with positive dimensional automorphism group.

\item Test the stability of $C_0$ against one-parameter subgroups $\rho$ coming from $\Aut(C_0)$. 

\item Compute the basin of attraction of (the Hilbert/Chow point of) $C_0$ and show that it contains $C$.
\end{enumerate}
For the  degenerate $C_0$ and $\rho$ coming from $\Aut(C_0)$, we can {\it explicitly compute} the Hilbert-Mumford index. Consider the following facts:
\begin{enumerate}
\item The $\rho$-action is nontrivial only on the rational subcurve  from which $\rho$ comes. So we can compute the exact weight of the $\rho$-action by using the Gr\"obner technique;
\item $\mu([C_0]_m, \rho)$ is a quadratic polynomial for $m$ larger than or equal to the regularity of $C_0$. 
\end{enumerate}
Hence by computing $\mu([C_0]_m, \rho)$ for $m$ up to the regularity of $C_0$ plus a couple, we can explicitly compute the index and  hence  determine {\it exactly for which $m$} $C_0$ is stable, semistable or unstable with respect to $\rho$. 

\subsection{Basin of attraction and deformation space with $\bG_m$-action}\label{S:basin}
How does one compute the basin of attraction of the  degenerate $C_0$? Assuming that our parameter space $X$ is smooth (in fact, one can show using deformation theory that the relevant Hilbert and Chow parameter spaces are smooth), 
we have the Bia{\l}ynicki-Birula decomposition $X = \amalg X_i$ with the following properties:
\begin{enumerate}
\item The sets $X_i^\rho$ of points in $X_i$ fixed under the $\rho$-action are precisely the connected components of the fixed point set $X^\rho$.
\item There is an affine bundle morphism $\varpi_i : X_i \to X_i^\rho$ such that 
$\varpi_i^{-1}(x^\star) = A_\rho(x^\star)$, $x^\star \in X_i^\rho$.
\item $T_{x^\star}X_i = (T_{x^\star}X)_{\ge 0}$, $
x^\star \in X_i^\rho$.
\end{enumerate}
In the third item, the subscript $\ge 0$ means the non-negative weight subspace of the $\bG_m$ action via $\rho$. The decomposition suggests an effective method to compute the basin of attraction. Indeed, let $X$ be the Hilbert scheme $\Hilb_{g,n}$  and $C_0$ be a curve with positive dimensional automorphism group. Under certain stability conditions, for instance c- or h-semistability \cite[Definition~2.5, 2.6]{HH2}, $\Hilb_{g,n}$ is smooth at $x^\star:=[C_0]$ and we apply the Bia{\l}ynicki-Birula decomposition theorem. To figure out which curves are in the basin of attraction $A_\rho(x^\star)$, we analyze the $\bG_m$ action on the tangent space $T_{[C_0]}\Hilb_{g,n}$. Actually for our purpose, by Luna's \'etale slice theorem, we may work with a $\bG_m$-invariant locally closed \'{e}tale slice $W \subset \Hilb_{g,n}$ containing $x^\star$. 
But in this case, the space $\Def(C_0)$ of first order deformations at $0$ is $\bG_m$-equivariantly isomorphic to $W$ at $x^\star$ \'etale locally, and the problem is now reduced to analyzing the $\bG_m$ action on the deformation space. The analysis in \cite{HM, HH2} all follow this line of ideas and basin of attraction computation is done simply by checking which curves are in the non-negative weight subspace of the deformation space. In fact, we only looked at the local deformation spaces at the singular points, in view of the exact sequence
\[
0 \to LT^1(C_0) \to T^1(C_0) \to \prod_{p\in \Sing(C_0)} T^1(\widehat{\cO}_{C_0, p}) \to 0
\]
where $T^1$ denotes the vector space of first order deformations and $LT^1(C_0)$ is the subspace of locally trivial deformations.
As seen in \cite{HH2}, for curves with at worst nodes, cusps and tacnodes as singularities, $\bG_m$ acts trivially on  $LT^1(C_0)$  and the problem was reduced further to analyzing the $\bG_m$ action on $\prod_{p\in \Sing(C_0)} T^1(\widehat{\cO}_{C_0, p})$. This is no longer the case in higher order singularities because of crimping \cite[Section~6]{ASvdW}.

As an application, we analyze the $m$-Hilbert stability of maximally degenerate curves with $A_{2b}$-singularities.  We will employ two techniques---a direct GIT computation of the Hilbert-Mumford indices and an alternative approach by explicitly computing the induced characters of the automorphism group on natural line bundles.

\subsection{Computing the Hilbert-Mumford index of cuspidal tails}\label{S:tails}
\begin{proposition}\label{P:tails}\cite{H}
 Let $C = D\cup_p R$ be a bicanonical genus $g$ curve such that $R$ is a rational curve of genus $b \ge 2$ with an $A_{2b}$ singularity  $y^2 = x^{2b+1}$ and meets $D$ in a node $p$. Then there exists a one parameter subgroup $\rho$ of $SL_{3g-3}$ coming from ${\rm Aut}(R)$ such that
\[
\mu([C]_m, \rho) =  \frac13(m-1) ((4b^2-8b+2)m - 3b^2)
\]
Retain $C$ and $\rho$, and consider the basin of attraction $A_\rho([C]_m)$.
We obtain:
\begin{enumerate}
\item Let $C' = D \cup_{p'} R'$ be a bicanonical curve where $R'$ is a hyperelliptic
curve of genus $b\ge 2$ meeting $D$ in a node $p'$, and $p'$ is a Weierstrass point of $R'$. Then there exists a 1-PS $\rho$ such that
\[
\mu([C']_m, \rho) = \frac13(m-1) ((4b^2-8b+2)m - 3b^2).
\]
 $R'$ (and sometimes $C'$ itself by abusing terminology) is called a
\emph{Weierstrass genus $b$ tail}.
\item Let $C''$ be a bicanonical genus $g$ curve obtained from $C$ by replacing $R$ by
an $A_{2b}$ singularity. That is, $C''$ is of genus $g$, has $A_{2b}$ singularity at $p''$ and admits a
partial normalization $\nu : (D, p) \to (C'', p'')$. Then there exists a 1-PS $\rho$ such that
\[
\mu([C'']_m, \rho) = -\frac13(m-1) ((4b^2-8b+2)m - 3b^2).
\]
\end{enumerate}
\end{proposition}
This in particular implies that even though $C'$ is Deligne-Mumford stable (assuming $D$ is so), it is $m$-Hilbert unstable for $m \le \frac{3b^2}{4b^2-8b+2}$. For example, a genus two Weierstrass tail is $m$-Hilbert unstable for $m \le 5$ and at best strictly semistable for $m = 6$. 
\begin{proof}
$R$ is a rational curve of genus $b$ with a single cusp $q$ whose local analytic equation is $y^2 = x^{2b+1}$.  $C = D\cup_p R$ is a bicanonical curve of genus $g$ consisting of $R$ and a genus $g-b$ curve $D$ meeting in a single node $p$. Restricting $\cO_C(1)$ to $R$ (resp. $D$), we find that it is of degree $4b-2$ (resp. $4g-4b-2$) and contained in a linear subspace of dimension $3b-1$ (resp. $3g-3b-1$). We can and shall choose coordinates such that
$R \subset \{x_{3b-1} = x_{3b} = \cdots = x_{3g-4} = 0\}$ and $D \subset \{x_0 = x_1 = \cdots = x_{3b-3} = 0\}$. $R$ may be parameterized by mapping $[s,t]$ to
\[
[s^{4b-2}, s^{4b-4}t^2, s^{4b-6}t^4, \cdots, s^{2b-2}t^{2b}, s^{2b-3}t^{2b+1}, s^{2b-4}t^{2b+2}, \cdots, t^{4b-2}]
\]
so that $R$ has a single cusp $ x_{b+1}^{2} = x_1^{2b+1}$ at $ q = [1, 0, \dots, 0]$, where we abused notation and let $x_1$ and $x_{b+1}$ denote their images in the completion of the local ring at $q$. Let $\rho$ denote the one-parameter subgroup with weights $$(0,2,4,6,\dots,2b, 2b+1, 2b+2, \dots, 4b-2, 4b-2, \dots, 4b-2).$$
The sum of these weights is $r := \binom{4b-1}2 - b^2 + (4b-2)(3g-3b-2)$.
We shall fix the $\rho$-weighted GLex order on the monomials.
\begin{lemma} The sum $w_{R,\rho}(m)$ of the weights of the degree $m$ monomials in $x_0, \dots, x_{3b-2}$ that are not in the initial ideal of $R$ is
\[
w_{R,\rho}(m) = (8b^2-8b+2)m^2 + (2b-1)m - b^2.
\]
\end{lemma}

\begin{proof} In general, weight computation of this sort can be accomplished by using Gr\"obner basis, but in this case there is a more elementary solution since $R$ admits a parameterization. Let $P_R(m) = (4m-1)(b-1) + 2m$, the Hilbert polynomial of $R$.  A monomial of degree $m$ pulls back to one of the following $m(4b-2)+1-b$ monomials
\[
s^{m(4b-2)-i}t^i, \quad i = 0, 2, 4, \dots, 2b, 2b+1, 2b+2, \dots, m(4b-2).
\]
If $\prod_{i\in I,|I|=m} x_i$ and $\prod_{i \in J,|J|=m} x_i$ pull back to the same monomial, then $$\prod_{i\in I,|I|=m} x_i - \prod_{i \in J,|J|=m} x_i$$ is in the initial ideal $in_\rho(I_R)$ of the ideal $I_R$ of $R$ with respect to the $\rho$-weighted GLex order. It follows that each $s^{m(4b-2)-i}t^i$ appears at most once
among the pullbacks of degree $m$ monomials not in the initial ideal $in_\rho(I_R)$. Since $m(4b-2)+1-b$ equals $P_R(m)$, it has to appear in the set exactly once. Therefore,
\[
\begin{array}{clll}
w_{R,\rho}(m) & = & \sum_{k=0}^b 2k + \sum_{k=2b+1}^{m(4b-2)}k \\
& = & (8b^2-8b+2)m^2 + (2b-1)m - b^2.
\end{array}
\]
\end{proof}

On the other hand, the contribution from $D$ to the total weight is
\[
\begin{array}{clllll}
w_{D,\rho}(m) & = & (4b-2)m \cdot h^0((\cO_C^{\otimes 2}|D)^{\otimes m} (-p))\\
& = &(4b-2)m((4m-1)(g-b-1)+2m-1)
\end{array}
\]
since $\rho$ acts on $D$ trivially with constant weight $4b-2$. Combining these weights and the average weight, we obtain the Hilbert-Mumford index of $C$: \small
\[
\begin{array}{lllllllllll}
\mu([C]_m, \rho) & = & \frac{mP(m)}{N+1}r - w_{R,\rho}(m) - w_{D,\rho}(m)
 \\
& = & \frac13m(4m-1)(\frac12(4b-1)(4b-2) - b^2 + (4b-2)(3g-3b-2)) \\
&& - ((8b^2-8b+2)m^2 + (2b-1)m - b^2) \\
&& - ((4b-2)m((4m-1)(g-b-1)+2m-1))\\
& = & \frac13\left((4b^2-8b+2)m^2 + (-7b^2+8b-2)m + 3b^2\right) \\
& = & \frac13(m-1) ((4b^2-8b+2)m - 3b^2)
\end{array}
\]\normalsize

Next, we analyze the basin of attraction. The local versal deformation space of $q$ is given by
\[
x_{2b+1}^2 = x_1^{2b+1} + c_{2b-1} x_1^{2b-1} + c_{2b-2} x_1^{2b-2} + \cdots +c_0.
\]
Since $\rho$ acts on $x_1$ and $x_{2b+1}$ with weights $2$ and $2b+1$ respectively, it acts on $c_i$ with positive weight $4b+2-2i$, $0 \le i \le 2b-1$. Hence the basin of attraction contains arbitrary smoothing of the cusp $q$. By considering the local stable reduction  \cite[\S~6.2.2]{Has00}, we can deduce that if $D\cup_p R$ is in the basin $A_{\rho}([C]_m)$, then $R$ must be hyperelliptic and $p$ is a Weierstrass point of $R$. Indeed, consider the isotrivial family $\mathcal C \to B=\spec k[[t]]$ whose general member is $\rho(t). (D\cup_p R)$ and the special member is $C$. Stable reduction of $\mathcal C$ yields $\mathcal C' \to B'$, $B' \to B$ a finite covering, whose general member is isomorphic to that of $\mathcal C$ and the special member is $D\cup_{p'} R'$ where $R'$ is  hyperelliptic and $p'$ is a Weierstrass point of $R'$: This is precisely the content of \cite[\S~6.2.2]{Has00}. By the separateness of $\overline{\mathcal M}_g$, it follows that  $R$ is hyperelliptic and $p$ is a Weierstrass point. On the other hand, $\rho$ acts with weight $-1$ on the local versal deformation of the node $p$ and the basin of attraction does not contain any smoothing of the node $p$. The assertion of the item (2) follows since $\rho^{-1}$ acts with the opposite weight, and the basin of attraction contains arbitrary smoothing of the node $p$ but no smoothing of the cusp $q$. 

 \begin{figure}[ht]
\centerline{\scalebox{0.6}{\psfig{figure=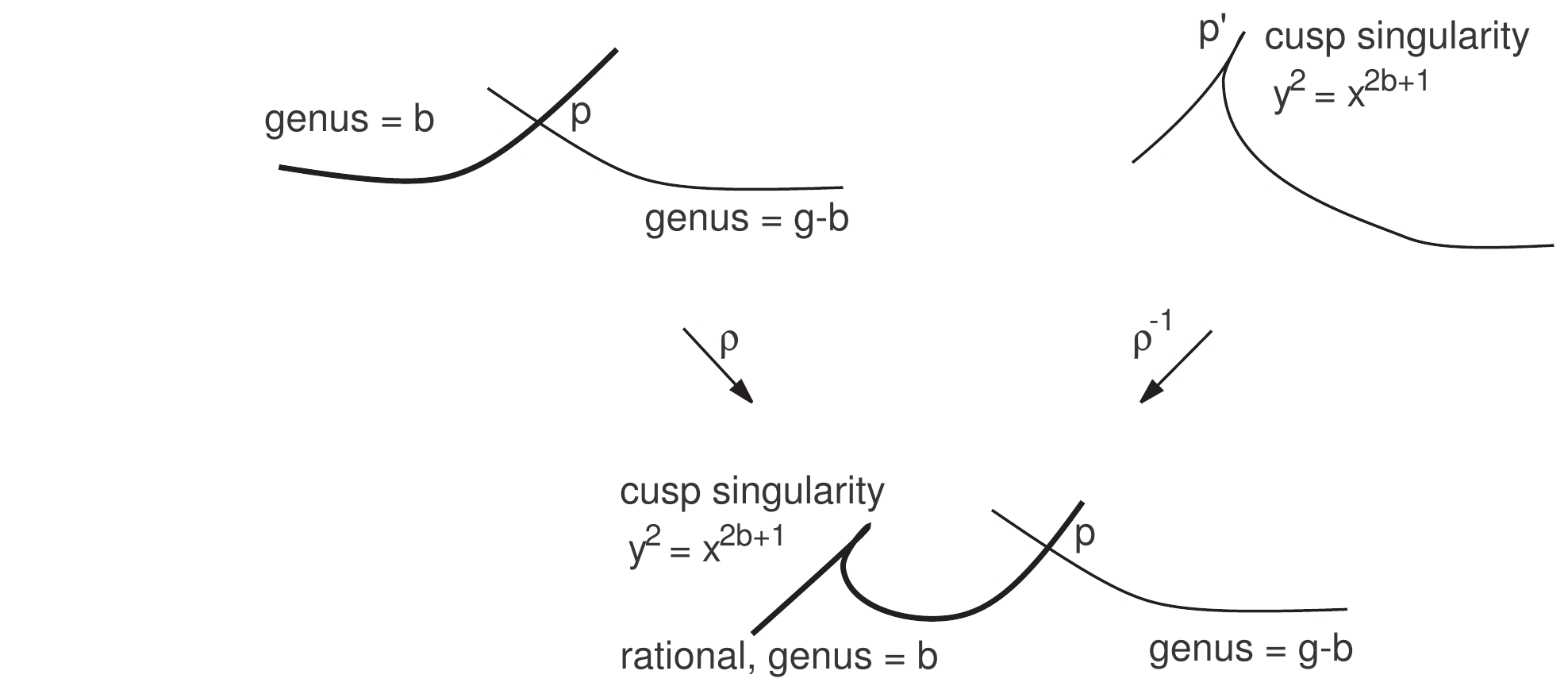}}}
\caption{Basin of attraction of cuspidal tail}
\label{F:cusps}
\end{figure}

This completes the proof of Proposition~\ref{P:tails}.
\end{proof}


 \subsection{An alternative approach using character theory}

As in Section \ref{S:tails}, let $C = D\cup_p R$ be a genus $g$ curve where $D$ is a smooth curve of genus $g-b$ and $R$, a rational curve meeting $D$ at a node $p$ with a monomial $A_{2b}$-singularity  $y^2 = x^{2b+1}$ at $q$.  One checks that there is an isomorphism $\eta: \mathbb{G}_m \to \Aut^{\circ}(C)$. Indeed, if we choose an isomorphism $\widetilde{R} \simeq \mathbb{P}^1$ such that 0 lies above the higher cusp and $\infty$ lies above the node, and let $t$ be a uniformizer at zero, then these automorphisms act on $R$ by $t \rightarrow \a t$ (and act trivially on $D$), where $\a$ is the coordinate of $\bG_m$.  As explained rigorously in \cite{AFS}, the line bundles $\lambda$, $\lambda_2$, $\delta$ and $K$ extend to a neighborhood of $[C]$ in the the stack of all curves and induce characters $\chi_{\lambda}(C, \eta)$, $\chi_{\lambda_2}(C, \eta)$, $\chi_{\delta}(C, \eta)$ and $\chi_{K}(C, \eta)$.  In this section, we will explicitly compute these characters and then recall how these characters gives an alternative computation of the Hilbert-Mumford indices computed in Section \ref{S:tails}.  (Of course, there are linear relations $K=13 \lambda -2 \delta$ and $\lambda_2 = 13 \lambda - \delta$ so that computing all four characters is redundant.)

To compute the character  $\chi_{\lambda}(C, \eta)$, consider a basis of $H^0(C, \omega_C)$ given by

$$
\left(0, \frac{dt}{t^{2b}}\right),\left(0, \frac{dt}{t^{2b-2}}\right), \ldots, \left(0, \frac{dt}{t^2}\right), (\omega_{1}, 0), \ldots, \left(\omega_{g-b},0 \right) ,
$$
where $\omega_{1}, \ldots, \omega_{g-b}$ is a basis for $\omega_{D}$.  Evidently $\mathbb{G}_m$ acts on this basis with weights $2b-1, 2b-3, \ldots, 3,  1, 0, \ldots, 0.$  Since $\lambda_1|_{[C]}=\bigwedge^{g} H^0(C, \omega_{C})$, we deduce that 
$$
\chi_{\lambda}(C, \eta)=\sum_{i=1}^{b}(2i-1)=b^2.
$$

Similarly, a basis for $H^0(C, \omega_C^2)$ is given by 
\begin{multline*}
\left(0, \frac{(dt)^2}{t^{4b}}\right),\left(0, \frac{(dt)^2}{t^{4b-2}}\right), \ldots, \left(0, \frac{(dt)^2}{t^{2b}}
\right), \left(0, \frac{(dt)^2}{t^{2b-1}}\right), \ldots, \left(\omega_{0}, \frac{(dt)^2}{t^{2}}\right),
\\ (\omega_{1}, 0), \ldots, \left(\omega_{3g-3b-2},0 \right) ,
\end{multline*}
where $\omega_{1}, \ldots, \omega_{3g-3b-2}$ is basis for $H^0(D,\omega_{D}^2(p))$, and 
$\omega_{0}$ is an appropriately chosen element of $H^0(D,\omega_{D}^2(2p)) \backslash H^0(D,\omega_{D}^2(p))$.  It follows that
$$\chi_{\lambda_2}(C, \eta)=\sum_{i=0}^{b-1}(2b+2i)+\sum_{i=0}^{2k-2}i =5b^2-4b+1.$$

To compute the $\chi_{\delta}(C, \eta)$ and $\chi_{K}(C,\eta)$, we write the first order deformation space
as
$$T^1(C) = T^1(D, p) \times \text{Cr}(\widehat{\mathcal{O}}_{C,q}) \times T^1(\widehat{\mathcal{O}}_{C,q}) \times
T^1(\widehat{\mathcal{O}}_{C,p})$$
where $\text{Cr}(\widehat{\mathcal{O}}_{C,q})$ denotes the ``crimping'' deformations (see \cite{ASvdW} for more details).  We can choose coordinates
$$\begin{aligned}
T^1(\widehat{\mathcal{O}}_{C,q}) &= \{ y^2 - x^{2b+1} + c_{2b-1}x^{2b-1} + \cdots + c_1 x +
c_0 = 0 \}, \\
T^1(\widehat{\mathcal{O}}_{C,p}) &= \{xy+ n=0\},
\end{aligned}$$
where $\mathbb{G}_m$
acts via $c_i \mapsto \lambda^{2i-4b-2} c_i$ and  $n \mapsto \lambda n$.  

By \cite[Proposition 5.7]{AFS}, the character $\chi_{\delta}(C,\eta)$ is simply the additive inverse of the weighted degree of the discriminant.  Since the discriminant of the $A_{2b}$-singularity has weighted degree $-4b(2b+1)$ while node has weighted degree $1$, we have
$$\chi_{\delta}(C,\eta) = 8b^2+4b-1.$$
By \cite[Lemma~4.15]{AFS}, the character $\chi_{K}(C,\eta)$ is the character of $T^1(C)$.   Using the above descriptions, we compute that the character of $T^1(\widehat{\mathcal{O}}_{C,q})$ is $-(4+6+\cdots + (4b+2)) = -(4b^2+6b)$.  The character of $T^1(\widehat{\mathcal{O}}_{C,p})$ is $1$.  The character of $ T^1(C_0, p)$ is trivial.  For $b \ge 2$, by \cite[Proposition 3.4]{ASvdW}, the weights of the action on $\text{Cr}(\widehat{\mathcal{O}}_{C,q})$ are $1,3, \ldots, 2b-3$.  Therefore, the character of $\text{Cr}(\widehat{\mathcal{O}}_{C,q})$ is $(b-1)^2$.  It follows that 
$$\chi_K(C, \eta) = -3b^2-8b+2.$$

We have therefore established:
\begin{proposition}  Let $C =  D\cup_p R$ be the genus $g$ curve with a nodally attached rational curve $R$ with a monomial $A_{2b}$-singularity.  Let $\eta: \mathbb{G}_m \to \Aut^{\circ}(C)$ be the isomorphism given above.  Then we have the following expressions for the characters:
$$\begin{aligned}
\chi_{\lambda}(C, \eta) &=b^2 \\
\chi_{\lambda_2}(C, \eta) &=5b^2-4b+1\\
\chi_{\delta}(C, \eta) &=8b^2+4b-1 \\
\chi_{K}(C, \eta) &= -3b^2-8b+2
\end{aligned}$$
It follows that the $m$th Hilbert-Mumford index is 
$$
\mu([C]_m, \rho) =  \frac13(m-1) ((4b^2-8b+2)m - 3b^2).
$$ 
\end{proposition}

\begin{proof} The final statement follows from the usual computation of the divisor class of the GIT polarization (see \cite[Proposition 7.1]{AFS}).
\end{proof}

\section{Local study of the moduli spaces of c-semistable and of h-semistable curves}\label{S:factorial}
In this section, we take a closer look at the flip \cite[Theorem~2.12]{HH2}
\[
\xymatrix{ \Mg(\frac7{10}+\e) \simeq \Mps \ar[dr]^-{\Psi} & &
\Mg(\frac7{10}-\e) \simeq \Mg^{hs}
\ar[dl]_-{\Psi^+}\\
& \Mg(\frac7{10})\simeq \Mg^{cs} &\\
}
\]
We will (1) give an \'etale-local description of this flip and (2) show that $\Mg^{hs}$ is not $\bQ$-factorial.

\subsection{ \'{E}tale local presentation of the flip} \label{S:etale-flip}
 An \'{e}tale local study of the Hassett-Keel program is rigorously carried out in \cite[Section~7]{ASvdW} where it is shown that \'etale  locally around any closed point $[C] \in
\bar{\cM}_g^{cs}$, the open inclusions of stacks
$$\bar{\cM}_{g}^{ps} \subseteq \bar{\cM}_{g}^{cs} \supseteq
\bar{\cM}_{g}^{hs}$$
correspond to the open chambers
$$\text{T}^1(C) ^- \subseteq \text{T}^1(C) \supseteq \text{T}^1(C) ^+$$ given by applying variation of GIT to the action of $\Aut(C)$ on the first order deformation space
$\text{T}^1(C)$.   We carry out this \'etale-local presentation in the particular case where $C$ is a general closed curve with infinite automorphism.  We then use this presentation to examine the flip in \cite{HH2}. 

 Let $C = D\cup R_1 \cup R_2$ be a c-semistable curve of genus $g$ consisting of 

\begin{enumerate}
\item a smooth curve $D$ of genus $g-2$;
\item smooth rational curves $R_i$ meeting each other in a tacnode $p$ and meeting $D$ in  a node $q_i$, $i = 1, 2$. 
\end{enumerate}

The first order deformation space of $C$ is
\[
\text{T}^1(C) \simeq \text{T}^1(D,q_1,q_2) \times \text{T}^1(\widehat\cO_{C,p}) \times \text{T}^1(\widehat\cO_{C,q_1}) \times \text{T}^1(\widehat{\cO}_{C,q_2})
\]
and there are isomorphisms
$$\begin{aligned}
\text{T}^1(\widehat{\mathcal{O}}_{C, p}) &= \{ y^2 = x^4 +  s_2x^2 + s_1x + s_0\,:\, s_i \in \mathbb{C} \} \\
\text{T}^1(\widehat{\mathcal{O}}_{C, q_i}) &= \{y^2 = x^2 + n_i\,:\, n_i \in \mathbb{C}\}.
\end{aligned}$$
If we fix an isomorphism $\Aut(C)^{\circ} \cong \mathbb{G}_m = \Spec
\mathbb{C}[\a]_\a$ which acts on a local coordinate $z$ around $0 \in \mathbb{P}^1$ by
$z \mapsto \a z$, then the action of $\Aut(C)^{\circ}$ on $\text{T}^1(C)$ is given
by
$$ s_i \mapsto \a^{-4+i}s_i, \qquad n_i \mapsto \lambda n_i$$
and is trivial on $\text{T}^1(D,p)$. Now, in this simple case, the open chambers $\text{T}^1(C)^-$ and $\text{T}^1(C)^+$ are defined as the non-vanishing locus of functions of negative and positive weight respectively. That is, $\text{T}^1(C)^- \subseteq \text{T}^1(C) \supseteq
\text{T}^1(C)^+$ are defined by the closed loci:
$$
\text{T}^1(C)  \setminus \text{T}^1(C)^- = V(s_0,s_1,s_2) \qquad \text{ and } \qquad
\text{T}^1(C)  \setminus \text{T}^1(C)^+  = V(n_1,n_2).
$$
The flip \'{e}tale-locally
 at $[C]$ corresponds to the flip arising from variation of GIT on the first order deformation space: 
\[
\xymatrix{\text{T}^1(C)^-/\!\!/ \bG_m  \ar[dr] & &
\text{T}^1(C)^+/\!\!/ \bG_m
\ar[dl]\\
& \Def(C)/\!\!/\bG_m &\\
}
\]
The locus $V(s_0,s_1,s_2)$ is the locus where $p$ remains a tacnode and $V(n_1,n_2)$ is the locus of curves containing an elliptic bridge.   This is precisely the picture obtained in \cite[Figure~1]{HH2} by carrying out the GIT. 
 \begin{figure}[ht]
\centerline{\scalebox{0.6}{\psfig{figure=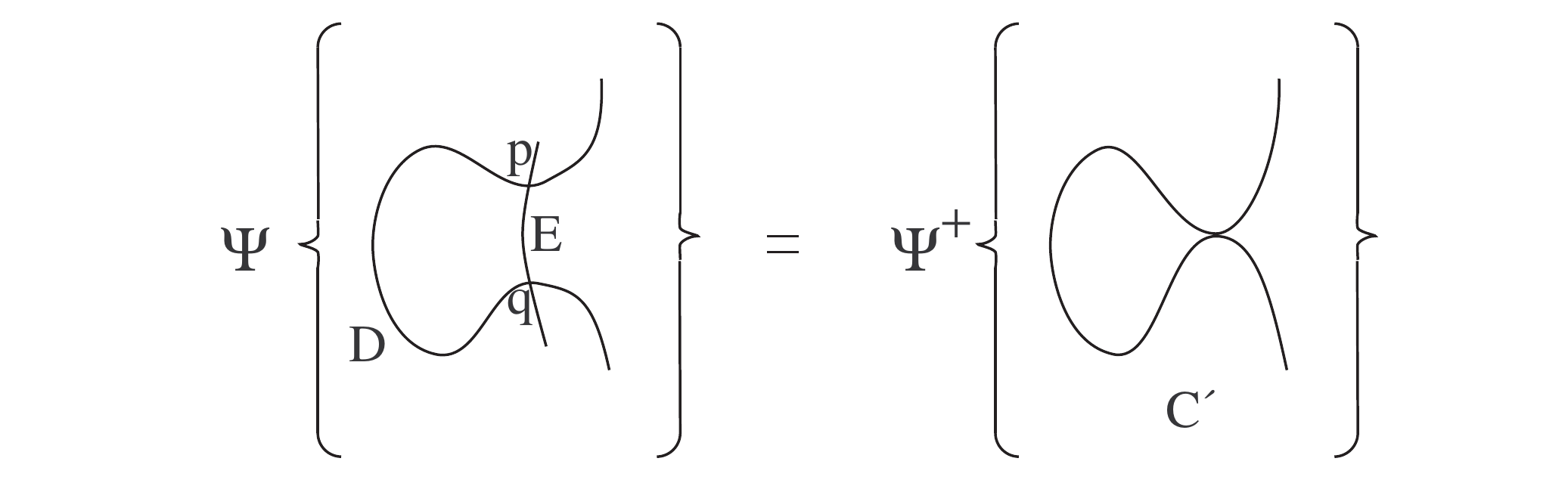}}}
\caption{The first flip}
\label{F:flip}
\end{figure}

\subsection{ $\bQ$-factoriality of log-canonical models}\label{sS:factorial}
We now show that $\Mg^{hs}$ is not $\bQ$-factorial for $g \ge 7$. This may seem unexpected since 
$\bQ$-factoriality is preserved under a flip in the log minimal model program as defined in \cite[Definition 2.8]{kollar-mori}.  However, the contraction 
$$\Psi: \Mg^{ps} = \bar{M}_g(7/10+\epsilon)  \to   \bar{M}_{g}(7/10)= \Mg^{cs}, $$
arising from the log minimal model program where the divisor $\delta$ is scaled, contracts more than one extremal ray.  Indeed, $\Psi$ contracts $\lfloor \frac{g-1}{2} \rfloor$ $(K+7/10\delta)$-negative extremal rays in $\bar{NE}(\bar{M}_g(7/10+\epsilon))$ corresponding to the ray $T_0$ of separating elliptic bridges and the rays $T_i$ for $i = 2, \ldots, \lfloor \frac{g-1}{2} \rfloor$ of elliptic bridges which separate a curve into a genus $i$ and $g-i-1$ components.  However, if one scales a generic boundary divisor $D = \sum_{i=0}^{\lfloor g/2 \rfloor} \alpha_i \delta_i$, then one expects $\lfloor \frac{g-1}{2} \rfloor$ flips corresponding to contracting each $B_i$ separately; see Remark \ref{remark-generic}.

To figure out the extremal rays contracted by $\Psi$, recall the divisorial contraction $T: \Mg \to \Mps$ which contracts the extremal ray $R$ generated by elliptic tails. The small contraction $\Psi: \Mg^{ps} \to \Mg^{cs}$ is induced by  $K_{\Mg^{ps}} + 7/10 \, \delta^{ps}$ and the ample cone of $\Mg^{ps}$ may be analyzed in terms of the ample cone of $\Mg$ via $T$.  We have
$$T^*\mathrm{NS}(\Mg^{ps})=R^{\perp} \subset \mathrm{NS}(\Mg),$$
the hyperplane spanned by a facet of the nef cone
of $\Mg$. 
$T^*(K_{\FMps} + 7/10\ \dps)$ contracts $\D_1$ and the following
one dimensional strata \cite{faber-1997, GKM}.  Let $X_0$ be a 4-pointed
stable curve of genus zero with one point moving and the other three
fixed.
\begin{itemize}
\item{ Attach two fixed 2-pointed curves of genus $1$ and
$g-3$, respectively, to $X_0$.  The extremal contraction of this ray has the locus $T_0$ of  elliptic bridges as the exceptional locus:
\[
T_0 = \{ C_1 \cup_{p,q} C_2 \, | \, g(C_1) = 1, g(C_2) = g-2\};
\]
}
\item{ Attach two 1-pointed curves of genus $i$ and $g-1-i$ (with $i \ge 1$ and $g-1-i \ge 1$)
respectively, and a 2-pointed curve of genus $0$ to $X_0$. This extremal contraction  has $T_i$ as
the exceptional locus: 
\[T_i = \{ C_1\cup_p C_2 \cup_q C_3 \, | \, g(C_1) =
i, g(C_2) = 1, g(C_3) = g-1-i \} .
\]
}
\end{itemize}
Since $\Psi$ contracts multiple extremal rays corresponding to various elliptic bridges, $\bQ$-factoriality may not necessarily be preserved under $\Psi$; indeed, by Theorem \ref{T:not-q-factorial}, $\bar{M}_g^{hs}$ is not $\bQ$-factorial for $g \ge 7$.

We will deduce the failure of $\bQ$-factoriality by applying the theorem below \cite{Drezet}.

\begin{theorem}\label{T:factorial} Let $Z$ be a nonsingular variety on which a reductive group $G$ acts, admitting a good quotient $\pi : Z \to M$. Suppose that there exists a saturated open subset $Z_0 \subset Z$ such that
\begin{enumerate}
\item $\codim_Z(Z\setminus Z_0) \ge 2$;
\item $\pi|_{Z_0} : Z_0 \to \pi(Z_0)$ is a geometric quotient;
\item $G_z$ is finite for $z \in Z_0$.
\end{enumerate}
Then $M$ is $\bQ$-factorial if and only if for every $G$-line bundle $L$ on $Z$ and every closed point $z \in Z$ with closed orbit $G\cdot z$, the stabilizer $G_z$ acts trivially on the fiber $L_z$. 
\end{theorem}

\begin{remark} \cite[Theorem~8.3]{Drezet} has stronger assumption that the action of $G$ on $Z_0$ be free so that $\pi(Z_0)$ is smooth and $\Pic(\pi(Z_0)) = \Cl(Z_0)$. Here, we are concerned about the $\bQ$-factoriality, and $\Pic(\pi(Z_0))_{\bQ} = \Cl(Z_0)_{\bQ}$ is valid with the finite stabilizer assumption.
\end{remark}

\begin{theorem} \label{T:not-q-factorial}
$\bar{M}_{g}^{hs}$ is $\mathbb{Q}$-factorial if and only if $g \le 6$.  
\end{theorem}

\begin{proof}
In our setting, we have $G = GL_{3g-3}(k)$ acting on $Z = \Hilb_{g,2}^{ss}$, the locus of Hilbert semistable curves.  By deformation theory, $Z$ is smooth.  Let $Z_0 = \Hilb_{g,2}^\circ \subseteq Z$ be the locus of Deligne-Mumford stable curves.  Then (1) $\codim_Z(Z\setminus Z_0) \ge 2$,  
(2) $\pi|_{Z_0}$ is a geometric quotient  and (3) $G_z$ is finite for $z \in Z_0$: (1) follows from  simple dimension calculation, and (2) and (3) follow since a Deligne-Mumford stable curve is Hilbert stable. 

For $g > 6$, Let $z  \in \Hilb_{g,2}/\!\!/G$ be the point representing the curve $C = C_1\cup R \cup C_2$, where
\begin{enumerate}
\item $
R = R_1\cup R_2\cup R_3$   is the union of rational curves where $R_1$ meets $R_2$ and $R_2$ 
meets $R_3$ at monomial $A_3$-singularities $p_2$ and $p_3$, respectively;
\item $ C_1$ (resp. $C_2$) is a
genus $2$ (resp. $g - 4$) curve meeting $R_1$ (resp. $R_3$) at a node.
\end{enumerate}
Then $z$ is a maximally degenerate curve corresponding to a closed orbit and $\Aut(C)^\circ \simeq  \bG_m$ acts on the fiber of $\delta_{2}$
with nontrivial character $\pm 1$. This follows from considering the $\bG_m$-action on the first order deformation space as in Section \ref{S:etale-flip}.  The divisor $\delta_2$ is cut out by the deformation parameter of the node (where $C_1$ meets $R_1$) which has weight $\pm 1$; it follows from  \cite[Proposition~5.7]{AFS} that the character of $\delta_2$ is $\pm 1$.

For $g < 7$, by an analysis of the maximally degenerate h-semistable curves and their induced characters, one can show that $\bar{M}_g^{hs}$ is $\mathbb{Q}$-factorial.
\end{proof}

 \begin{remark}  {\bf (Log canonical models arising from generic scalings)}  \label{remark-generic}
 
 Instead of studying the ``democratic" log minimal models $\bar{M}_g(\alpha)$, one can consider a generic log canonical model
 $$\bar{M}_g(D) :=
  \Proj \bigoplus_{d \ge 0} \Gamma(\bar{\cM}_g, \lfloor d( K_{\bar{\cM}_g} +D)\rfloor), \qquad \text{ where }D = \sum_{i=0}^{\lfloor g/2 \rfloor} \alpha_i \delta_i$$
for rational numbers $0 \le \alpha_i \le 1$.  The divisor $K + D$ is F-nef if $-9+12\alpha_0-\alpha_1 \ge 0$ and $\alpha_i + \alpha_j + \alpha_k + \alpha_l - \alpha_{i+j} - \alpha_{i+k} - \alpha_{i+l} \le 2$ for integers $i,j,k,l$ with $g=i+j+k+l$.  Moreover, by \cite[Prop. 6.1]{GKM}, $K+D$ is nef if $\alpha_i \le \alpha_0$ for $i \ge 1$.

Choose $D= \sum_{i=0}^{\lfloor g/2 \rfloor} \alpha_i \delta_i$ with  $\alpha_i \le \alpha_0$ for $i \ge 1$ such that $K + D$ is ample.  Then $K+tD$ is ample for $1 \ge t > 9/(12\alpha_0-\alpha_1)$.  One can use the character computations in \cite{AFS} to make predictions on the moduli interpretations of the log canonical models $\bar{M}_g(tD)$ arising from scaling $t$ from $1$ to $0$.   At $t = 9/(12\alpha_0-\alpha_1)$, one expects cusps to replace elliptic tails in $\bar{M}_g(tD)$.  At $t = t_0=7/(10 \alpha_0)$, one expects tacnodes to replace non-separating elliptic bridges.  For each $i = 2, \ldots, \lfloor \frac{g-1}{2} \rfloor$, at $t = t_i=7/(12 \alpha_0-\alpha_i - \alpha_{i+1})+ \epsilon$ (where in the case that $g=2k+1$ is odd, we use the convention that $\alpha_{k}=\alpha_{k+1}$ so that $t_k=7/(12 \alpha_0-2\alpha_k) + \epsilon$), one expects tacnodes to replace elliptic bridges which separate the curve into genus $i$ and $g-i-1$ components.  If $t_i \neq t_j$ for $i \neq j \in \{0, 2, 3, \ldots \lfloor \frac{g-1}{2} \rfloor\}$, one expects flips
$$\xymatrix{
\bar{M}_g((t_i+\epsilon)D) \ar[rd]	&	& \bar{M}_{g}((t_i-\epsilon)D)  \ar[ld] \\
&   \bar{M}_{g}(t_iD)
}$$
for $i=0,2,\ldots \lfloor \frac{g-1}{2} \rfloor$ with $\bar{M}_{g}((t_i-\epsilon)D)$ $\bQ$-factorial.
\end{remark}

\clearpage

\begin{table}[p]\label{T:1}
\begin{center}
\renewcommand{\arraystretch}{2}
\begin{tabular}{|p{1in}|p{1.in}|p{1in}|p{1in}|p{1in}|p{1in}|}
\hline 
   & Curve embedding & Linearization & Log canonical model & Stability $\&$ Singularity \\
\hline
 $\Hilb$& $n$-canonical, $n\ge 5$ & $m \gg 0$ & $\Mg(1) \simeq \Mg$ & Deligne-Mumford; $A_1$ \\
 \hline
$\Chow$ & four-canonical &  & $\Mg(9/11)$ & pseudo-stability; $A_1, A_2$\\
\hline

$\Hilb$ & four-canonical & $m\gg 0$  & $\Mg(9/11-\epsilon)$ & pseudo-stability; $A_1, A_2$\\
\hline

$\Chow$ & tri-canonical &   & $\Mg(25/32)$ & pseudo-stability; $A_1, A_2$\\
\hline

$\Hilb$ & tri-canonical & $m\gg 0$  & $\Mg(25/32-\epsilon)$ & pseudo-stability; $A_1, A_2$\\
\hline

$\Chow$ & bi-canonical &   & $\Mg(7/10)$ & c-semistability $A_1\sim A_3$\\
\hline

$\Hilb$ & bi-canonical & $m\gg 0$  & $\Mg(7/10-\epsilon)$ & h-semistability $A_1 \sim A_3$\\
\hline

$\Hilb$ & bi-canonical & $m = 6$  & $\Mg(2/3)$ &  $A_1 \sim A_4$\\
\hline

$\Hilb$ & bi-canonical & $m = 4.5$  & $\Mg(19/29)$ &  $A_1 \sim A_4, A_5^\dagger$\\
\hline

$\Hilb$ & bi-canonical & $m = 2.25$  & $\Mg(17/28)$ &  $A_1 \sim A_5$\\
\hline



$\Chow$ & canonical &   & $\Mg(\frac{3g+8}{8g+4})$ & ADE, ribbons, etc\\
\hline


\end{tabular}

\caption{Log canonical models as GIT quotients \cite{FS}}
\end{center}
\end{table} 

\clearpage

\bibliographystyle{alpha}
\bibliography{survey}
\end{document}